\documentclass[11pt]{amsart}

\usepackage{amssymb}
\usepackage{epsfig}
\usepackage{comment}
\usepackage{amsmath}
\usepackage{subfigure}
\usepackage[section]{placeins}
\usepackage{mathrsfs}
\usepackage{graphicx}
\usepackage[notrig]{physics}
\usepackage{units}
\usepackage{pstricks}
\usepackage{color}
\usepackage{algorithm}
\usepackage{algorithmic}
\usepackage{multirow}
\definecolor{orange}{rgb}{1.0, 0.549, 0}

\newtheorem{thm}{Theorem}[section]

\theoremstyle{definition}

\theoremstyle{remark}
\newtheorem{remark}[thm]{Remark}
\graphicspath{{./}{./Figures}}

\begin{document}

\title[The finite-particle method]{An optimal-transport finite-particle method for driven mass diffusion}

\author
{
A.~Pandolfi$^1$, I.~Romero$^{2,3}$ and M.~Ortiz$^4$
}

\address
{
${}^1$Politecnico di Milano,
Civil and Environmental Engineering Department, 20133 Milano, Italy
\\
${}^2$Universidad Polit\'ecnica de Madrid,
Mechanical Engineering Department, 28006 Madrid, Spain
\\
${}^3$IMDEA Materials Institute,
28906 Madrid, Spain
\\
${}^4$California Institute of Technology,
Engineering and Applied Science Division, Pasadena CA, 91125, USA
}

\email{anna.pandolfi@polimi.it, ignacio.romero@upm.es, ortiz@caltech.edu}

\begin{abstract}
We formulate a finite-particle method of mass transport that accounts for general mixed boundary conditions. The particle method couples a geometrically-exact treatment of advection; Wasserstein gradient-flow dynamics; and a Kullback-Leibler representation of the entropy. General boundary conditions are enforced by introducing an adsorption/depletion layer at the boundary wherein particles are added or removed as dictated by the boundary conditions. We demonstrate the range and scope of the method through a number of examples of application, including absorption of particles into a sphere and flow through pipes of square and circular cross section, with and without occlusions. In all cases, the solution is observed to converge weakly, or in the sense of local averages.
\end{abstract}

\maketitle


\section{Introduction}
\label{sec:Introduction}

In a previous publication \cite{Pandolfi:2023}, we have developed a class of velocity-free finite-particle methods by recourse to a time-discretized variational formulation of mass diffusion \cite{JordanKinderlehrerOtto1997, JordanKinderlehrerOtto1998, JordanKinderlehrerOtto1999, LiHabbalOrtiz2010, fedeli2017geometrically}, and by discretizing the mass density by means of finite particles, or blobs, as in \cite{carrillo2017, carrillo2019}. In this optimal-transport framework, the motion of the particles is the result of a competition between entropy, as given by the Kullback-Leibler functional, and mobility, expressed in terms of the Wasserstein distance (cf., e.~g., \cite{Villani2003}) between consecutive configurations. Thus, entropy works to disperse the particles uniformly over the domain, whereas mobility works to hinder their motion. In addition, a fractional-step decomposition allows advective terms and sources to be integrated exactly. This geometrically-exact integration confers the approach great robustness, in contrast with methods based on upwinding or stabilization, which remain challenging to date.

The work presented in \cite{Pandolfi:2023} focuses on forging a clear and direct link between the theory of transport of measures and particle methods as approximations thereof, while the implementation of general boundary conditions is given scant attention. In the present work, we extend the implementation of the finite-particle method of \cite{Pandolfi:2023} to general mixed boundary conditions. The extension accounts for equilibrium conditions with the environment through an external chemical potential and, therefore, it enables consideration of physically-based Dirichlet boundary conditions. The extension also accounts for mass exchange with the environment through the enforcement of Neumann inlet/outlet boundary conditions. These two types of particle exchange with the environment are represented by means of a further fractional step decomposition consisting of a predictor step, unconstrained by the boundary conditions, followed by a corrector step, which brings the solution back into balance with the boundary conditions. The all-important corrector step is formulated by introducing a boundary layer wherein particles are added or removed in accordance with the boundary conditions.

We showcase the range and scope of the method through a number of examples of application, including absorption of particles into a sphere and flow through pipes of square and circular cross section, with and without occlusions. In all these examples, the numerical solutions are observed to converge weakly to the exact solution in the sense of measures.

\section{Physical basis of mass transport in solids}

In order to frame the general class of problems under consideration and make contact with applications, we succinctly summarize the thermodynamical basis of mass transport in solids.

\subsection{General thermodynamic framework}

Let $\rho \geq 0$ denote the mass density of a certain diffusive species. Thermodynamics postulates the existence of a {\sl Gibbs free energy density} per unit volume $G(\rho, T)$, where $T$ is the absolute temperature that is minimized at equilibrium \cite{callen1985wt}. The chemical potential $\mu$ is conjugate to the concentration and defined as
\begin{equation} \label{eq-chemical-pot}
    \mu = \frac{\partial G}{\partial\rho} (\rho, T) .
\end{equation}
In the absence of sources mass balance requires
\begin{equation} \label{eq-bmass}
    \frac{\partial\rho}{\partial t} + \nabla\cdot J = 0,
\end{equation}
where $J$ is the mass flux and $\nabla\cdot$ denotes the divergence operator.

Consider now an arbitrary subdomain $U$. Under isothermal conditions, the power due to the supply of mass to $U$ is
\begin{equation} \label{eq-expended}
    \mathcal{P}
    =
    -
    \int_{\partial U} \mu \, J \cdot \nu \, dA ,
\end{equation}
where $\nu$ is the outward unit normal to the boundary $\partial U$ of $U$ and ${J}\cdot \nu$ is the outgoing mass flux. The rate of dissipation in $U$ is, therefore,
\begin{equation} \label{eq-dissipation}
    \mathcal{D}
    =
    \mathcal{P}
    -
    \frac{d}{dt} \int_U G(\rho, T) \, dV ,
\end{equation}
where $T$ is held constant. An appeal to the chain rule and the divergence theorem gives
\begin{equation} \label{eq-dissipation2}
    \mathcal{D}
    =
    -
    \int_U
    \Big(
        \nabla \cdot (\mu J )
        +
        \frac{\partial G}{\partial \rho} \,
        \frac{\partial\rho}{\partial t}  \Big)
    \, dV ,
\end{equation}
which by (\ref{eq-chemical-pot}) and (\ref{eq-bmass}) reduces to
\begin{equation} \label{eq-dissipation3}
    \mathcal{D}
    =
    -
    \int_U
        J \cdot \nabla \mu
    \, dV .
\end{equation}
Since $U$ is arbitrary, it follows from this identity that the dissipation is non-negative, as required by the second law, if and only if the dissipation inequality
\begin{equation} \label{eq-mf-dissip}
    J \cdot\nabla\mu \leq 0
\end{equation}
is satisfied. We note from this inequality that the driving force for the mass flux is the gradient of the chemical potential. Suppose next that there is a kinetic potential $\varphi(\rho, T; \nabla\mu)$ such that
\begin{equation}
    J
    =
    -
    \frac{\partial\varphi}{\partial\nabla\mu}(\rho, T; \nabla\mu) .
\end{equation}
Then, the dissipation inequality (\ref{eq-mf-dissip}) follows if $\varphi(\rho, T; \cdot)$ is a convex function and $\varphi(\rho, T; 0) = 0$.

\subsection{Ideal solid solutions}
\label{Subsec:Diffusion}

By way of example, consider an impurity/host system in the form of a dilute ideal solution in which the impurity phase diffuses in a pressure and potential energy field. Under these conditions, the free energy density per unit volume of the solid is \cite{serebrinsky2004vi, dileo2013fu}:
\begin{equation}
    G(\rho, T)
    =
    \rho \, \mu_0(T)
    +
    R T \,
    \Big(
        \rho \log\frac{\rho}{\rho_{\rm ref}} - \rho
    \Big)
    +
    \rho ( U - \, p V ) ,
\end{equation}
where $\rho$ is the mass density of the impurity phase, $\mu_0(T)$ is a temperature-dependent reference chemical potential, $p(x)$ is the pressure field, $U(x)$ is a potential energy per unit mass, $R$ is the universal gas constant, $\rho_{\rm ref}$ is a reference mass density, $V$ is the partial molar volume of the impurity in solid solution, and we omit the dependencies on position for economy of notation. Using Eq.~\eqref{eq-chemical-pot}, the corresponding chemical potential evaluates to
\begin{equation} \label{Eq:BulkMu}
    \mu
    =
    \mu_0(T)
    +
    R T\,
    \log\frac{\rho}{\rho_{\rm ref}}
    -
    p V
    +
    U .
\end{equation}
Assuming a dissipation potential of the form
\begin{equation}
    \varphi(\rho, T; \nabla\mu) = \frac{\rho M(T)}{2} \, |\nabla\mu|^2,
\end{equation}
where $M(T)$ is a temperature-dependent mobility, we obtain Fick's law
\begin{equation} \label{Eq:Fick1}
    J = - \rho \, M(T)  \nabla \mu .
\end{equation}
Inserting (\ref{Eq:BulkMu}) into (\ref{Eq:Fick1}) we obtain
\begin{equation} \label{Eq:Fick1b}
    J =
    -
    \kappa(T) \, \nabla \rho
    +
    \rho u(T)
    -
    \rho \nabla \Psi(T),
\end{equation}
where
\begin{equation}
    \kappa(T) = M(T) \, R T ,
    \quad
    u(T) = \frac{\kappa(T) V}{R T} \nabla p ,
    \quad
    \Psi(T) = \frac{U}{M(T)}
\end{equation}
are the bulk diffusivity, Darcy drag velocity due to the pressure gradient, and the drag velocity potential due to the potential energy field, respectively. Inserting (\ref{Eq:Fick1b}) into \eqref{eq-bmass} we finally obtain the Fokker-Planck equation
\begin{equation} \label{Eq:Fick2}
    \frac{\partial \rho}{\partial t}
    +
    \nabla
    \cdot
    (\rho u)
    =
    \kappa \, \Delta \rho
    +
    \nabla \cdot (\rho \nabla \Psi) ,
\end{equation}
which jointly accounts for advection, diffusion, and energetic forces.

\section{Problem definition}
\label{sec:ProblemDefinition}

The preceding physical considerations lead to the general advection-diffusion initial-boundary-value mass-transport problem
\begin{subequations}\label{system_complete}
\begin{align}
    &
    \partial_t \rho + \nabla\cdot(\rho u)
    =
    \kappa \, \Delta \rho
    +
    \rho s,
    & \text{in } \Omega \times [0,T] ,
    \label{eq:TD:Diff1}
    \\ &
    \rho
    =
    g ,
    & \text{on } \Gamma_D \times [0,T] ,
    \label{eq:TD:Diff2}
    \\ &
    \kappa \nabla\rho \cdot \nu
    =
    f ,
    & \text{on } \Gamma_N \times [0,T] ,
    \label{eq:TD:Diff3}
    \\ &
    \rho(x,0) = \rho_0(x),
    & \text{in } \Omega ,
    \label{eq:TD:Diff4}
\end{align}
\end{subequations}
where $\Omega \subset \mathbb{R}^d$ is the domain of analysis and $[0,T]$ the time interval under consideration, $\nu$ is the outward unit normal to the boundary $\Gamma$, $\rho$ is the unknown mass density, $\rho_0$ is its initial value and $\kappa$ is the mass diffusivity.
In writing (\ref{system_complete}) and henceforth, we assume isothermal conditions and omit all reference to temperature.

In addition to the bulk equilibrium and kinetics relations discussed in the foregoing, in formulating problem (\ref{system_complete}), we have also introduced advection by a given transport velocity field $u$, and distributed mass sources per unit volume $\rho s$. As accounted for in (\ref{eq:TD:Diff3}), the advection velocity contributes to the mass flux through the Neumann boundary. However, for thermodynamic consistency and to avoid ill-conditioning from over-prescribed boundary conditions, we assume
\begin{equation}
    u \cdot \nu = 0, \quad \text{on } \Gamma , 
\end{equation}
i.~e., the advection velocity field is tangent to the boundary and defines a flow that maps the domain $\Omega$ to itself.

We specifically consider mixed boundary conditions of two forms: thermodynamic equilibrium with an external chemical potential, Eq.~(\ref{eq:TD:Diff2}); and prescribed outward mass flux, Eq.~(\ref{eq:TD:Diff3}). We recall that the condition of thermodynamic equilibrium at the boundary is $\mu = \mu_{\rm ext}$, where $\mu_{\rm ext}$ is the external chemical potential and $\mu$ obeys the equation of state (\ref{eq-chemical-pot}). Assuming the dependence of $\mu$ on $\rho$ to be one-to-one, we obtain (\ref{eq:TD:Diff2}) with $g = \mu^{-1}(\mu_{\rm ext})$ as the equilibrium density at the boundary.

We note that the resulting boundary condition (\ref{eq:TD:Diff2}) is of the Dirichlet type, and $\Gamma_D$ may thus be regarded as the Dirichlet boundary. The remaining boundary condition (\ref{eq:TD:Diff3}), which enforces a prescribed outward flux $f$ through the boundary, is of the Neumann type, and $\Gamma_N$ may thus be regarded as the Neumann boundary. For well-posedness, we assume that $\Gamma_D \cap \Gamma_N = \emptyset$ and $\Gamma_D \cup \Gamma_N = \Gamma$.

Eq.~(\ref{eq:TD:Diff1}) can be equivalently reformulated as the transport problem
\begin{subequations}\label{eq:TD:RV}
\begin{align}
    & \label{eq:TD:RV1}
    \partial_t \rho + \nabla\cdot(\rho v)
    =
    \rho s ,
    \\ & \label{rVre5S}
    \rho v
    =
    \rho u
    -
    \kappa \nabla\rho ,
\end{align}
\end{subequations}
where $v$ is a velocity field that results from the combined effect of advection and diffusion, Eq.~(\ref{eq:TD:RV1}) is the continuity equation of Eulerian continuum mechanics and (\ref{rVre5S}) is a mobility law combining the effects of advection and diffusion. This reformulation clearly evinces the transport character of the problem and motivates its framing within a measure-theoretical---rather than functional space---framework.

\section{Approximation}
\label{sec:approximation}

We aim to approximate the solutions of (\ref{system_complete}) variationally through a combination of time discretization followed by spatial discretization of the mass density into particles. To this end, we exploit the additive structure of (\ref{eq:TD:Diff1}) to decompose the incremental updates using the method of fractional steps \cite{fedeli2017geometrically}. The advective fractional step can be conveniently solved exactly following the approach described in \cite{fedeli2017geometrically, Pandolfi:2023}. This exact geometrical treatment of advection confers particle methods great robustness. The source fractional can also be solved exactly pointwise \cite{fedeli2017geometrically, Pandolfi:2023}.

We may therefore focus attention to the diffusive fractional step, obtained formally by setting $u = 0$ and $s=0$ in \eqref{eq:TD:RV1}, with the result
\begin{subequations}
\begin{align}
    &
    \partial_t \rho
    =
    \kappa \, \Delta \rho ,
    & \text{in } \Omega  ,
    \\ & \label{Sfd3dY}
    \rho
    =
    g ,
    & \text{on } \Gamma_D ,
    \\ & \label{mV6Lpz}
    \kappa \nabla\rho \cdot \nu
    =
    f ,
    & \text{on } \Gamma_N .
\end{align}
\end{subequations}
The overall strategy for the approximation of solutions is to regard the boundary conditions (\ref{Sfd3dY}) and (\ref{mV6Lpz}) as {\sl constraints} and solve the constrained evolution problem by a projection method.

\subsection{Unconstrained step}

Let $[t_k,t_{k+1}]$ be a generic time step. We start the update by dropping the constraints (\ref{Sfd3dY}) and (\ref{mV6Lpz}) and solving an unconstrained diffusion equation over the entire space $\mathbb{R}^d$. The particles are held in the domain $\Omega$ by the addition of a penalty potential
\begin{equation}\label{Pe18A9}
    \Psi(x) = \frac{C}{2} {\rm dist}^2(x,\Omega) ,
\end{equation}
where
\begin{equation}
    {\rm dist}(x,\Omega)
    =
    \min\{ \|x-y\| \, : \, y \in \Omega \}
\end{equation}
is the distance function to $\Omega$ and $C$ is a penalty stiffness to be calibrated.

We characterize the attendant evolution of the density by means of the energy-dissipation functional \cite{JordanKinderlehrerOtto1997, JordanKinderlehrerOtto1998, JordanKinderlehrerOtto1999, fedeli2017geometrically}
\begin{equation}\label{eq:TD:F}
 \begin{split}
    &
    F(\rho_{k+1})
    =
    \frac{1}{2}
    \frac{D^2(\rho_k,\rho_{k+1})}{t_{k+1}-t_k}
    + 
    \int_{\mathbb{R}^d}
        \Big[
            \Psi
            +
            \kappa
            \log\Big(\frac{\rho_{k+1}}{\rho_{\rm ref}}\Big)
        \Big]
    \,\rho_{k+1}\, dx
    \to \inf!,
\end{split}
\end{equation}
where $D(\rho_k,\rho_{k+1})$ is the Wasserstein distance between $\rho_k$ and $\rho_{k+1}$. It is readily verified that the stationary points of (\ref{eq:TD:F}) satisfy a time-discretized form of the diffusion problem (\ref{system_complete}), e.~g., \cite{fedeli2017geometrically}.


It would seem natural to further approximate $\rho(x,t)$ by concentrating mass on a finite collection of points $\{x_p(t)\}_{p=1}^n$, i.~e., by setting
\begin{equation}\label{8EYqqq}
    \rho(x,t) \sim \sum_{p=1}^n m_p \delta(x-x_p(t)) ,
\end{equation}
where $\delta$ is the Dirac delta and $m_p$ is the mass carried by particle $p$. However, the functional (\ref{eq:TD:F}) is not defined for trial densities of this type. In \cite{fedeli2017geometrically}, this difficulty is circumvented by replacing the incremental functional (\ref{eq:TD:F}) with an alternative weak statement, modeled after (\ref{eq:TD:RV}), in which the mass density appears linearly and can therefore be approximated by the {\sl ansatz} (\ref{8EYqqq}).

In the present work, following Carrillo {\sl et al.} \cite{carrillo2017, carrillo2019} we instead spread, or `fatten', the particles over a finite width to be determined, for the purpose of evaluating the functional (\ref{eq:TD:F}). The corresponding approximation is
\begin{equation}\label{bC4KRV}
    \rho(x,t)
    \sim
    \sum_{p=1}^n m_p  \varphi_p(x - x_p(t), t) ,
\end{equation}
where the function $\varphi_p(x - x_p(t); t)$ represents the time-dependent profile of particle $x_p$. As in \cite{carrillo2017,carrillo2019}, we specifically choose a Gaussian particle profile of the form
\begin{equation}\label{eq:nfunction}
    \varphi_p(x - x_p(t), t)
    =
    \Big( \frac{\beta_p(t)}{\pi} \Big)^{d/2}
    \, {\rm e}^{-\beta_p(t) \| x - x_p(t) \|^2}
    :=
    {N}(x - x_p(t), \beta_p(t)),
\end{equation}
normalized so that
\begin{equation}
    \int_{\mathbb{R}^d} {N}(x - x_p(t), \beta_p(t)) \, dx = 1 .
\end{equation}
We note that ${N}(x - x_p(t), \beta_p(t))$ has dimensions of an inverse volume, carries unit mass and has a width set by the parameter $\beta_p(t)$.

\begin{remark}[$\rho$ is a measure]
It seems tempting to regard (\ref{bC4KRV}) as an interpolation scheme and the resulting $\rho(x,t)$ as a function in some linear space of functions. However, this interpretation is misleading and should be carefully avoided in practice. Instead, the object $\rho(x,t)$ resulting from (\ref{bC4KRV}) is to be regarded as a regular measure resulting from a mollification of the singular measure (\ref{8EYqqq}). Operationally, this means that the only operation that is natural on $\rho(x,t)$ is taking expectations, or averages, of continuous functions, i.~e.,
\begin{equation} \label{RG5QKm}
    E_\rho[f;t] = \int f(x) \, \rho(x,t) dx ,
\end{equation}
where $f(x)$ is a continuous function. In particular, the pointwise evaluation of $\rho(x,t)$ is essentially meaningless. Furthermore, evaluations of averages such as (\ref{RG5QKm}) can only be expected to be accurate when $f(x)$ is slowly varying on the scale of separation between particles.
\hfill$\square$
\end{remark}

\begin{remark}[Regularization]
Likewise, it seems tempting to regard the spatial discretization as a Galerkin scheme in which representation (\ref{bC4KRV}) is inserted into the functional (\ref{eq:TD:F}). Again, this interpretation is misleading and should be avoided. The correct interpretation is that (\ref{bC4KRV}) induces a regularization $F_\beta$ of the functional $F$. In particular, the value of $F_\beta$ at the discrete measure (\ref{8EYqqq}) equals the value of $F$ at the mollified measure (\ref{bC4KRV}). Convergence must therefore be understood in the sense of density of the discrete measures (\ref{8EYqqq}) as the number $N$ of particles and, simultaneously, in the sense of $\Gamma$-convergence of the regularized functionals $F_\beta$ to $F$ as $\beta \to +\infty$. Conditions on $N$ and $\beta$, as well as error estimates that establish the convergence of the approximation, can be found in \cite{Conti:2023}.
\hfill$\square$
\end{remark}


The Gaussian profile of the mollified particles and the absence of mass redistribution renders the calculation of the dissipation functional in (\ref{eq:TD:F}) straightforward \cite{Dowson:1982}. The choice of Gaussian profiles also facilitates the computation of integrals by numerical quadrature \cite{hildebrand2013introduction}. Combining these preceding approximations, the reduced energy-dissipation functional takes the form
\begin{equation}\label{ZXTj1Y}
\begin{split}
    F(\rho_{k+1})
    & \sim
    \sum_{p=1}^n
    \frac{1}{2}
    \frac
    {\| x_{p,k+1} - x_{p,k} \|^2 + ( \beta_{p,k+1}^{-1/2} - \beta_{p,k}^{-1/2} )^2}{t_{k+1}-t_k}
    \, m_p
    \\ & +
    \sum_{p=1}^n
    \kappa
    \log
    \Big(
        \frac{\rho_{k+1}(x_{p,k+1})}{\rho_{\rm ref}}
    \Big)
    \, m_p
    + 
    \sum_{p=1}^n
    \Psi(x_{p,k+1})\, m_p
    \to \min!
\end{split}
\end{equation}
with
\begin{equation}\label{2IeJn6}
    \rho_{k+1}(x)
    :=
    \sum_{p=1}^n
    \Big( \frac{\beta_{p,k+1}}{\pi} \Big)^{d/2}
    \,{\rm e}^{- \beta_{p,k+1} \| x - x_{p, k+1} \|^2 }
    \, m_p .
\end{equation}
We note that the incremental variational principle (\ref{ZXTj1Y}) determines both the updated position of the particles and their optimal widths. The corresponding Euler-Lagrange equations are
\begin{equation}\label{9eBxCb}
    \frac{\partial F}{\partial x_{r,k+1}}
    =
    0 ,
    \qquad
    \frac{\partial F}{\partial \beta_{r,k+1}}
    =
    0 ,
\end{equation}
which evaluate to
\begin{subequations}\label{j4K0Ht}
\begin{align}
\begin{split}
    & \label{bXq5IE}
    m_r
    \frac{x_{r,k+1} - x_{r,k}}{t_{k+1}-t_k}
    + \\ &
    \frac{\partial}{\partial x_{r,k+1}}
    \sum_{p=1}^n
    \left[
        \kappa
        \log
        \Big(
            \frac{\rho_{k+1}(x_{p,k+1})}{\rho_{\rm ref}}
        \Big)
        +
        \Psi(x_{p,k+1})
    \right]
    \, m_p
    =
    0 ,
\end{split}
    \\ & \label{W0gjO2}
    -
    \frac{m_r}{2}
    \frac
    {
        \beta_{r,k+1}^{-1/2}-\beta_{r,k}^{-1/2}
    }
    {
        \beta_{r,k+1}^{3/2} (t_{k+1} - t_k)
    }
    +
    \frac{\partial}{\partial \beta_{r,k+1}}
    \sum_{p=1}^n
        \kappa
        \log
        \Big(
            \frac{\rho_{k+1}(x_{p,k+1})}{\rho_{\rm ref}}
        \Big)
    \, m_p
    =
    0 .
\end{align}
\end{subequations}
These equations set forth an implicit time-stepping scheme for the positions and widths of the particles. In the following, $m_p$ is assumed constant for all the particles.

\subsection{Projection step}\label{sYFG37}

Evidently, the result $x_{k+1}^{\rm pre}$ of the preceding diffusive step may violate the constraints (\ref{Sfd3dY}) and (\ref{mV6Lpz}) in general. We reinstate these constraints by projecting the unconstrained predictor onto the corresponding constraint sets, carefully avoiding any artificial clustering effects that might arise from precluding mass flow at the boundary.

To this end, we partition the boundary $\Gamma$ into subsets $\mathcal{A} = \{A_l,\ l=1, \dots, L\}$ in a manner consistent with the partition of $\Gamma$ into the Dirichlet and Neumann boundaries $\Gamma_D$ and $\Gamma_N$, respectively. In addition, we define an boundary layer $\Gamma \times [-b,b]$ of thickness $2b$ and partition it into the subsets $\mathcal{B} = \{B_l,\ l=1, \dots, L\}$, with
\begin{equation}
    B_l
    =
    \{
        x\in\Gamma \times [-b,b] \, : \,
        A_l \text{ closest to }  x
    \} .
\end{equation}
Thus, the partition $\mathcal{B}$ of the boundary layer $\Gamma \times [-b,b]$ consists of cylindrical sets transversal to the boundary of $\Omega$ with base in $\mathcal{A}$. Each set $B_l$ can be split into an internal an external component defined respectively as
\begin{equation} \label{eq-bl-split}
    B_l^{int} = B_l\cap \Omega = A_l\times[-b,0]
    \ ,\qquad
    B_l^{ext} = B_l\setminus B_l^{int} = A_l\times(0,b]\ .
\end{equation}

Over the Dirichlet boundary $\Gamma_D$, the target mass in each domain $B_l^{int}$ at time $t_{k+1}$ is
\begin{equation}
    m_{l,k+1}
    \sim
    b \,
    \int_{A_l}
        g(x,t_{k+1})
    \, d\mathcal{H}^{d-1}(x) ,
\end{equation}
where $d\mathcal{H}^{1}(x)$ denotes the element of length and $d\mathcal{H}^{2}(x)$ denotes the element of area. The actual mass $m_{l,k+1}^{\rm pre}$  at the end of the diffusive step is likely to differ from $m_{l,k+1}$. This discrepancy is eliminated simply by inserting $\lfloor(m_{l,k+1}-m_{l,k+1}^{\rm pre})/m_p\rfloor$ particles if $m_{l,k+1} > m_{l,k+1}^{\rm pre}$, e.~g., at random locations within $B_l^{int}$; or by deleting $\lfloor(m_{l,k+1}^{\rm pre}-m_{l,k+1})/m_p\rfloor$ particles if $m_{l,k+1}^{\rm pre} > m_{l,k+1}$, e.~g., again randomly within $B_l^{int}$.

Driven by diffusive forces, the particles in $B_l^{int}$ can move toward the interior of $\Omega$ or crawl into $B_l^{ext}$. If the barrier potential $\Psi$ is extended to the outer boundary of $B_l^{ext}$, i.~e., on $\Gamma\times\{b\}$, then the domain $B_l^{ext}$ will effectively behave as a mass repository or buffer for the Dirichlet boundary condition at $A_l$, isolating outside $\Omega$ the spurious layer of particles that inevitably cluster near the impermeable barrier effected by the potential~$\Psi$.

The Neumann boundary conditions demand a different treatment. Over the Neumann boundary $\Gamma_N$ the target mass in each domain $B_l$ at time $t_{k+1}$ is
\begin{equation}
\label{eq-neu}
    m_{l,k+1}
    \sim
    m_{l,k}
    +
    \int_{t_k}^{t_{k+1}}
    \int_{A_l}
        f(x,t)
    \, d\mathcal{H}^{d-1}(x) \, dt .
\end{equation}
Consequently, the mass $m_{l,k+1}-m_{l,k}$ has to be inserted at every time step in each $B_l$. As before, this can be accomplished by inserting the appropriate number of particles at random positions of the domain, or by randomly eliminating enough particles so as to satisfy~(\ref{eq-neu}).

\subsection{Summary of the algorithm}

\begin{algorithm}
\begin{algorithmic}[1]
\REQUIRE Assign domain $\Omega$, initial density $\rho_0$, particle mass $m_p$, and target number of particles at the steady state $n_\infty$
\REQUIRE Estimate the particle distance
$\Delta r \sim \left( |\Omega|/n_\infty \right)^{1/d}$
\REQUIRE Assign width $\beta \sim 1/\Delta r^2$ and boundary layer thickness $b \sim \sqrt{L \Delta r}$
\REQUIRE Fill $\Omega$ with $n_0$ particles at random positions $x_{p,0}$ to obtain $\rho_0$
\REQUIRE Discretize the total time $T$ in $K$ steps
$\Delta t_k = t_{k+1} - t_k \leq \Delta r^2/\kappa$
\REQUIRE Approximate the mass densities as $\rho_0(x_{p,0})$, Eq.~\eqref{2IeJn6}
\FOR {$t = t_k$, $k=1,\ldots,K$}
\STATE Advection step: update the particle positions according to $u$
$$x^u_{p,k+1} = x_{p,k} + u(t_k) \Delta t_k$$
\STATE Source step: adjust the number of particles to account for sources
$$\rho_{k+1} = \rho_{k} \exp { s(t_k) \Delta t_k}$$
\STATE Unconstrained diffusion step: update the particle position
$$x_{r,k+1} = x^u_{r,k+1} +
\Delta t_k
\left[
\frac{\kappa \nabla \rho_k (x_{r,k})}{\rho_k(x_{r,k})}
+
\sum_{p=1}^n
\frac{\partial \Psi(x_{p,k})}{\partial x_{r,k}} \right] $$
\FOR {$l=1,\ldots, L$}
\IF {\rm Dirichlet boundary}
\STATE Compute the particle number $n^{D,l}$ giving density $g_l(x,t_{k+1})$
$$n^{D,l} = \lfloor m^{\rm pre}_{l,k+1}/m_p \rfloor$$
\STATE Count the number of particles $n^l$ in $B_l^{int}$
\IF {$n^l < n^{D,l}$}
\STATE Introduce $n^{D,l} - n^l$ particles in $B_l^{int}$
\ELSE
\STATE Eliminate $n^l - n^{D,l}$ random particles from $B_l^{int}$
\ENDIF
\ENDIF
\IF {Neumann boundary}
\STATE Estimate the number of inlet/outlet particles $n^{N,l}$
$$n^{N,l} = \left\lfloor (m_{l,k+1}-m_{l,k})/m_p \right\rfloor$$
\IF {$n^{N,l} > 0$}
\STATE Introduce $n^{N,l}$ in $B_l$
\ELSE
\STATE Eliminate $n^{N,l}$ from $B_l$
\ENDIF
\ENDIF
\ENDFOR
\ENDFOR
\end{algorithmic}
\caption{Fractional Step Diffusion Algorithm}
\label{alg:Diffusion}
\end{algorithm}

We note that considerable simplification is achieved if the width of the particles is assumed to be uniform, i.~e., $\beta_{r,k+1} = \beta_{k+1}$, and if the particle widths are assumed to be steady state at all times, i.~e., $\beta_{k+1} = \beta =$ constant. We expect convergence to be dependent on a careful choice of $\beta$ as a function of the number of particles~\cite{Conti:2023}.

The accuracy of the approximation depends critically on a careful choice of the parameters $\beta,b, \Delta t$ as a function of the number of particles, spatial dimensions, domain geometry, and material parameters.

A rigorous error analysis \cite{Conti:2023} shows that the optimal choice of $\beta$ is
\begin{equation}
  \label{eq-betadr}
    \beta \sim \Delta r^{-2} ,
\end{equation}
where
\begin{equation}
    \Delta r \sim {\left(\frac{|\Omega|}{n}\right)}^{1/d}
\end{equation}
is a length of the order of the mean separation between particles. The error, in the sense of the flat norm, in the approximation of a mass density $\rho$ by finite particles, or `blobs', is then of order $\beta^{-1/2}$, which proves a universal approximation property of finite particles.

In calculations, we further choose the thickness of the boundary layer to be intermediate between $\Delta r$ and the size of the specimen $L$. Specifically, we choose
\begin{equation}
    b \sim \sqrt{L\,\Delta r} .
\end{equation}

In order to evaluate the boundary conditions accurately, we choose the time step according to the Courant-Friedrichs-Lewy condition \cite{courant1928}
\begin{equation}
    \label{eq-stable-dt}
    \Delta t\le \Delta t_{\rm stable} \approx \frac{\Delta r^2}{\kappa} ,
\end{equation}
which is commensurate with the diffusive time scale. Finally, we choose the penalty stiffness $C$ in the barrier potential (\ref{Pe18A9}) to be the highest possible that does not decrease the stable time step, namely,
\begin{equation}
    \label{eq-penalty}
    C = \frac{1}{\Delta t} ,
\end{equation}
which concludes the setup of the model.

\section{Numerical examples}

Next, we present selected examples of application that demonstrate the range and scope of the finite-particle method and its convergence properties. For all the considered cases, the parameter $\beta$, the time step size $\Delta t$, and the penalty constant $C$ are chosen as indicated in Eqs.~(\ref{eq-betadr}), (\ref{eq-stable-dt}), and~(\ref{eq-penalty}), respectively.

\subsection{Dirichlet problem for a spherical domain}
\label{ssec:sphere}

\begin{figure}[ht]
    \centering \includegraphics[width=0.7\textwidth]{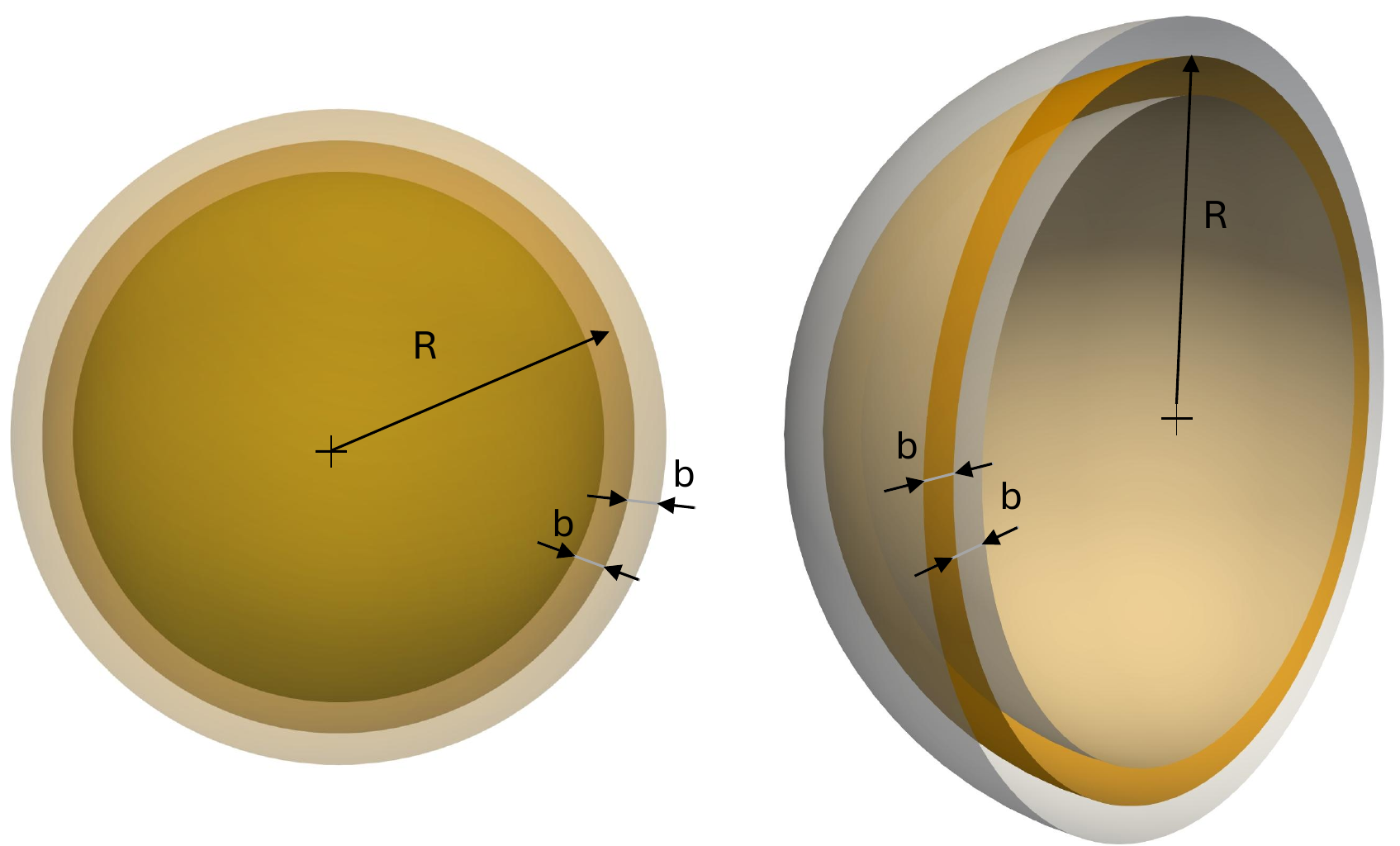}
    \caption{\footnotesize Geometry of the spherical domain used for the verification of the method. The image shows the double boundary layer used to enforce the Dirichlet boundary conditions. The thickness of the boundary is assigned as $2b$, of which $b$ inside the sphere and $b$ outside the sphere.}
    \label{fig-sphere-bc}
  \end{figure}

\begin{figure}[ht]
\includegraphics[width=0.7\textwidth]{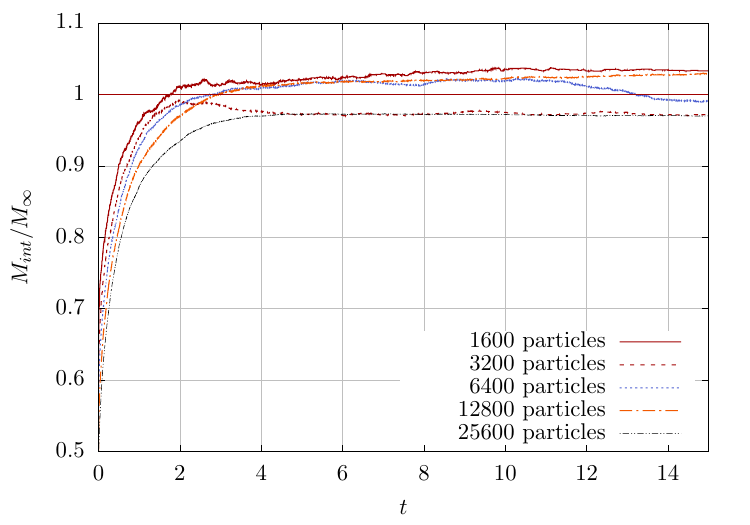}    \caption{\footnotesize
      Sphere example. Time evolution of the relative total mass inside the sphere for various numbers of particles~$n_\infty$.}
    \label{fig-sphere-m}
\end{figure}

\begin{figure}[ht] \includegraphics[width=0.7\textwidth]{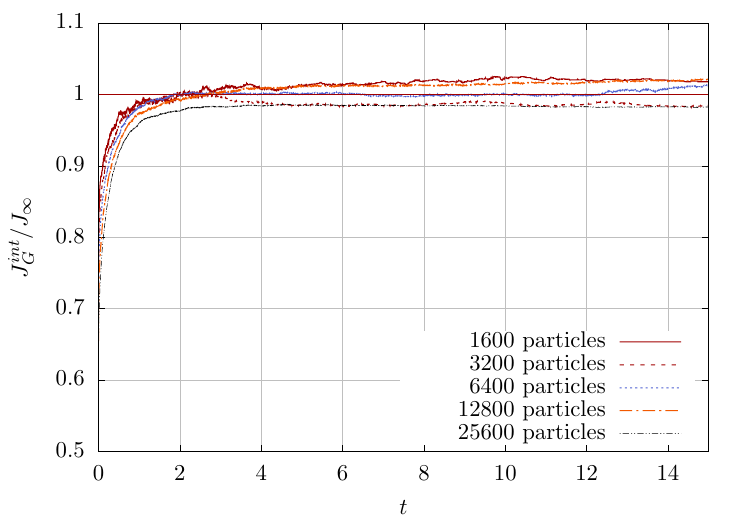}
    \caption{\footnotesize
      Sphere example. Time evolution of the relative polar inertia of the sphere for various numbers of particles~$n_\infty$.}
    \label{fig-sphere-j}
\end{figure}

In the first example, we consider the transport of mass within a spherical domain of radius~$R=1$ initially empty (see Fig.~\ref{fig-sphere-bc}). Uniform Dirichlet boundary conditions
are applied on the surface of the sphere, inducing a mass influx that progressively fills the domain until the density,
throughout the whole sphere, becomes homogeneous and equals the boundary value. For this simple geometry and choice of boundary conditions, the surface set is $\mathcal{A}\equiv\Gamma$ and  only one boundary layer $\mathcal{B}$ must be considered.
We set the diffusivity $\kappa=1$, the total mass at steady state $M_\infty= 1$, and the total analysis time $T=15$. The steady-state density $\rho_\infty$ follows as
\begin{equation}
  \label{eq-rhoinf}
    \rho_\infty = \frac{3 M_\infty}{4 \pi R^3} = 0.238732\, .
\end{equation}
We impose the boundary density $g$ to be equal to the reference density
$\rho_\infty$. Thus, the boundary condition will
enforce the final mass inside the sphere to be~$M_{\infty}$.

We choose a set of problems with an increasing number of particles at steady state, i.~e.,  $n_\infty=2^h\cdot 1600$ with $h=0,1,\ldots,4$, and in each case we assign the mass of one particle to be
\begin{equation}
    m_p = \frac{M_\infty}{n_\infty}\, .
\end{equation}

As suggested in Section~\ref{sec:approximation}, for the Gaussian particle profile we adopt $\beta=2/\Delta r^2$, with $\Delta r = R\,n_{\infty}^{-1/3}$. The particles are injected or removed with the help of a boundary layer of thickness $2b = 2\sqrt{R\,\Delta r}$. A spherical shell of thickness $b$ inside the spherical domain is used to gauge at each instant the density in the domain of analysis; another spherical shell of the same thickness, but external to the domain, is introduced with a barrier potential at radius $R+b$, see Fig.~\ref{fig-sphere-bc}.

For the chosen particle numbers, Fig.~\ref{fig-sphere-m} shows the time evolution of the mass $M(t)$ inside the sphere, defined as
\begin{equation}
\label{eq-mt}
 M(t) = \sum_{\substack{p=1\\ r_p(t)\le R}}^n m_p \, ,
\end{equation}
where $r_p(t)$ is the distance, at time $t$, of particle $p$ to the center of the sphere, excluding the contributions from the external spherical shell of
thickness~$b$. Fig.~\ref{fig-sphere-m} shows that, after an initial transient during which the sphere is filled, the total mass of the sphere approaches the theoretical value~$M_{\infty}$.

\begin{figure}[ht]
    \centering
    \includegraphics[width=0.7\textwidth]{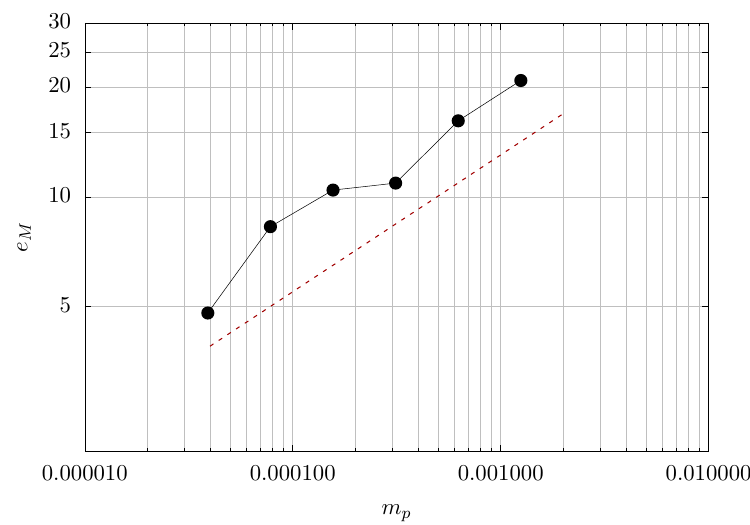}
    \caption{\footnotesize Sphere example.
    Error in the $L^1$ norm of the total mass~\eqref{eq-l1-sphere1}. The broken line represents the best-fitting power law showing convergence (error $e_{M}\approx m_p^{0.38}$).}
    \label{fig-sphere-l1m}
\end{figure}

\begin{figure}[ht]
    \centering
    \includegraphics[width=0.7\textwidth]{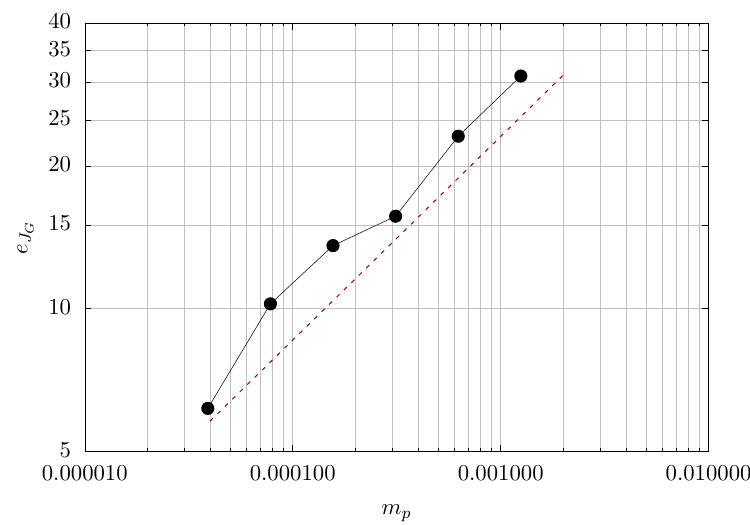}
    \caption{\footnotesize Sphere example.
    Error in the $L^1$ norm of the total polar inertia~\eqref{eq-l1-sphere2}. The broken line represents the best-fitting power law showing convergence (error $e_{J_G}\approx m_p^{0.43}$).}
    \label{fig-sphere-l1j}
\end{figure}

\begin{figure}[!t]
    \centering
    \includegraphics[width=0.8\textwidth]{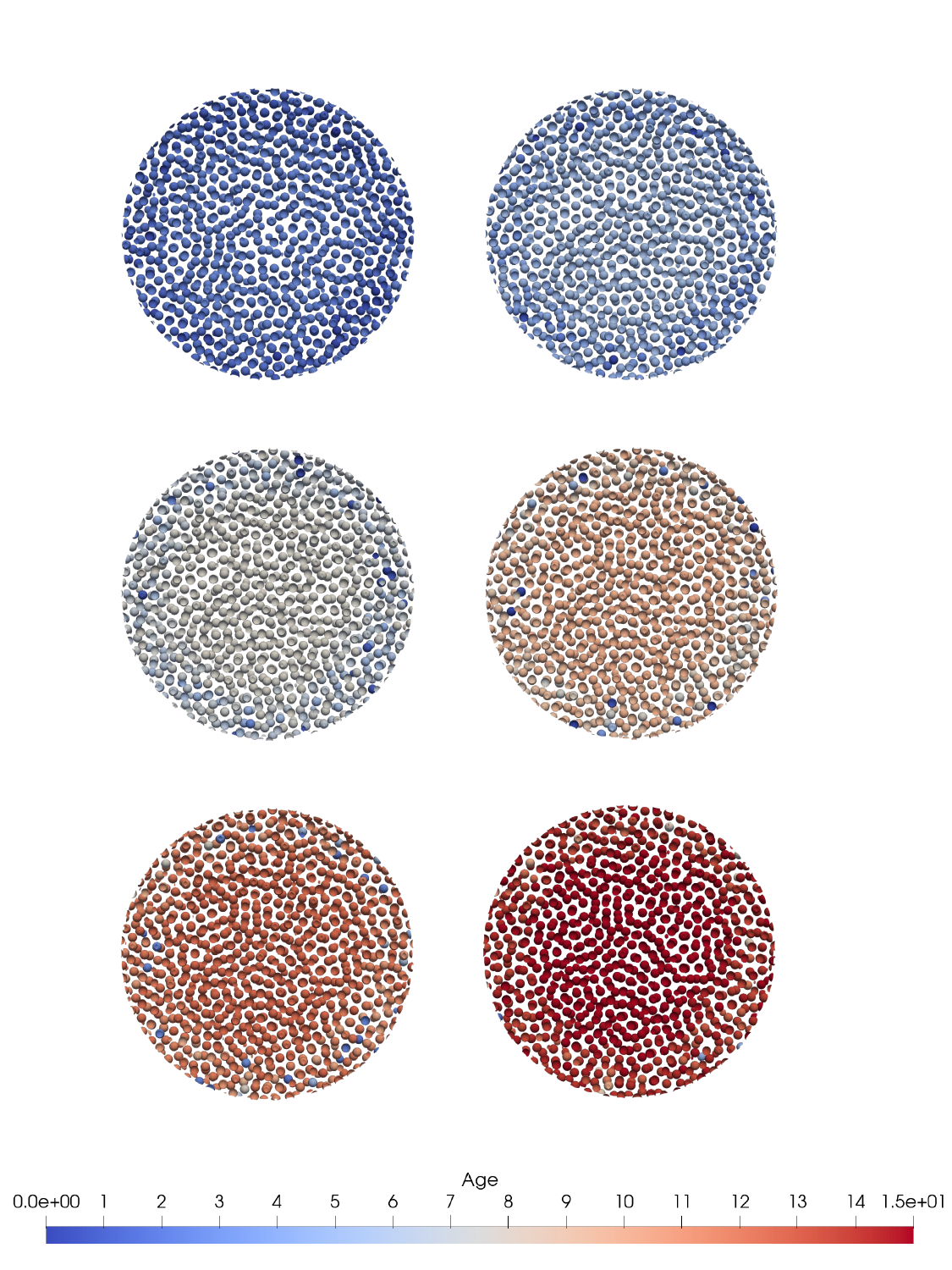}  \caption{\footnotesize Sphere example.
      Snapshots of a slice of thickness $b$ passing through the center of the sphere and clipping the particle distribution
      during the diffusion process ($n_p=6400$). From top to bottom, left to right, $t=2.5, 5.0, 7.5, 10.0, 12.5, 15.0$. Colors refer to particle age since insertion.}
    \label{fig-sphere-snap}
\end{figure}

To assess the uniformity of the particle distribution in the domain, we calculate the time evolution of the polar inertia with respect to the sphere center of the particles residing inside the sphere, defined as
\begin{equation}
    J_G(t) = \sum_{\substack{p=1\\ r_p(t)\le R}}^n r_p^2(t)\, m_p \, .
\end{equation}
As in the case of the mass, calculations ignore the contribution from the particles in the external
boundary layer. Fig.~\ref{fig-sphere-j} provides an illustration of the evolution of $J_G(t)$ relative to the polar inertia $J_\infty$ of a homogeneous sphere of density $\rho_{\infty}$
\begin{equation}
    J_\infty = \frac{3}{5} M_{\infty} R^2 = \frac{4}{5} \pi R^5 \rho_{\infty}\, .
\end{equation}
As for the total mass, after an initial transient the inertia approaches the value~$J_{\infty}$.
We note that for both the mass and the inertia calculations the limit values fluctuate and do not coincide with their theoretical values. This results from the explicit nature of the time integrator and the enforcement of the boundary conditions. The barrier potential prevents particles affected by diffusive forces from leaving the domain, while also introducing a short-range inhomogeneity in the mass density near the boundary. Placing the potential barrier at a distance $R+b$ from the center of the sphere reduces this effect, while minimizing the additional computational cost associated with tracking and accounting for interactions with the extra particles in the outer spherical shell. Remarkably, a thicker outer shell would reduce the spurious effect at the expense of a higher computational
cost. For the sphere example and all the cases considered, 

\begin{table}[h]
    \centering
    \begin{tabular}{rcccccrr}
        \hline
         $n_\infty$ & $m_p$ & $\Delta r$ & $\beta$ & $n_{\rm t}$ & $n_{\rm t}/n_{\infty}$ \\
         \hline
      \hline
       1600 & $   6.25\cdot10^{-4}$ & 0.085498797 & 273.59615 &  4130 &   2.58 \\
       3200 & $  3.125\cdot10^{-4}$ & 0.067860440 & 434.30682 &  6806 &   2.13 \\
       6400 & $ 1.5625\cdot10^{-4}$ & 0.053860867 & 689.41911 & 12688 &   1.98 \\
      12800 & $ 7.8125\cdot10^{-5}$ & 0.042749399 & 1094.3846 & 25154 &   1.97 \\
      25600 & $3.90625\cdot10^{-5}$ & 0.033930220 & 1737.2273 & 43138 &   1.69 \\
         \hline
    \end{tabular}
    \vskip 5pt
    \caption{\footnotesize
    Numerical parameters and results for the sphere example with Dirichlet boundary. $n_{\rm t}$ is the total
    number of particles, including the external layer, at the end of the simulation, $n_{\rm t}/n_{\infty}$ is the
    particle overhead due to the Dirichlet condition on the outer boundary.}
    \label{tab-sphere}
\end{table}

Table~\ref{tab-sphere} collects parameters and results relevant to the stationary state of the sphere with Dirichlet boundary, including the total number $n_t$ of particles used in numerical simulations. Note that the relative overhead $n_t/n_\infty$ decreases as the number of particles increases, since the ratio of the external shell's volume to the total volume approaches zero as $\Delta r \to 0$.

To verify the convergence of the method, we analyze the limit values of well-defined quantities of the mass density, such as the $L^1$ norms (in time) of mass and polar inertia, i.~e.,
\begin{subequations}
  \label{eq-l1norms-sphere}
  \begin{align}
    L^1_{M} &= \int_0^T M(t)\, dt \approx
              \sum_{k=0}^K \sum_{p=1}^{n(t_k)} m_p\,\Delta t\ ,
              \label{eq-l1-sphere1}
    \\
    L^1_{J_G} &= \int_0^T J_G(t)\, dt \approx
                \sum_{k=0}^K \sum_{p=1}^{n(t_k)} r_p^2(t)\, m_p \,\Delta t\ .
                \label{eq-l1-sphere2}
   \end{align}
\end{subequations}
Note that in \eqref{eq-l1norms-sphere} the sums must include all the mass points, including those in the external shells. Assuming that both quantities converge as $m_p\to0$, and that their corresponding errors can be expressed as power laws of $m_p$, the limiting values and convergence exponents can be estimated. The results of these computations are summarized in Figs.~\ref{fig-sphere-l1m}-\ref{fig-sphere-l1j}, which confirm the convergence norm of the mass and polar inertia in the $L^1$, respectively.

Fig.\ref{fig-sphere-snap} illustrates the mass diffusion process by showing snapshots of the particle distribution inside the sphere. The particles are colored according to their age, defined as the difference between the snapshot time and the corresponding insertion time. Older particles move toward the center of the sphere, while relatively younger ones accumulate near the boundary layer, where they may be inserted or removed, as explained in Section~\ref{sec:approximation}.

\subsection{Mass storage in a box with mixed boundary conditions}
\label{ssec:massStorage}

\begin{figure}[h]
    \centering
    \includegraphics[width=0.8\textwidth]{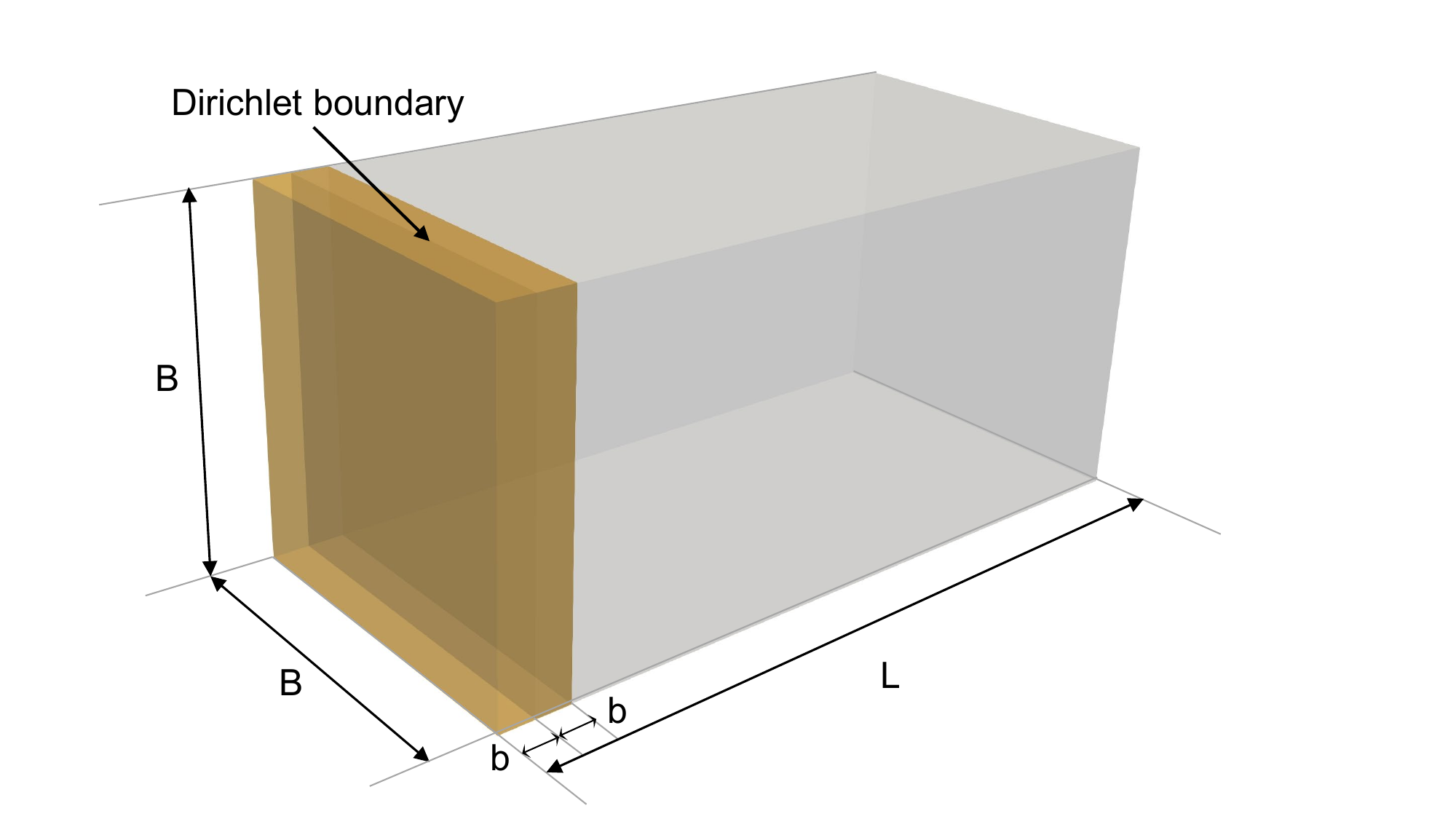}
    \caption{\footnotesize
      Parallelepiped problem. Dirichlet boundary conditions are applied to the left side of the parallelepiped.}
    \label{fig:fig7}
\end{figure}

Next, we simulate the filling of a container of parallelepiped shape. The box, initially empty, has two equal dimensions, $B=1$, and a larger dimension, $L=2$, see Fig.~\ref{fig:fig7}. Mass flux is allowed across one face, where the density is prescribed at a constant value $\rho_\infty=500$, while the barrier potential~\eqref{Pe18A9} is employed for the five no-flux faces. Following a transient phase, the container reaches a uniform density, with total mass $M_\infty = \rho_\infty |\Omega| = \rho_\infty B^2 L = 1000$. 

\begin{figure}[h]
    \centering
    \includegraphics[width=0.7\textwidth]{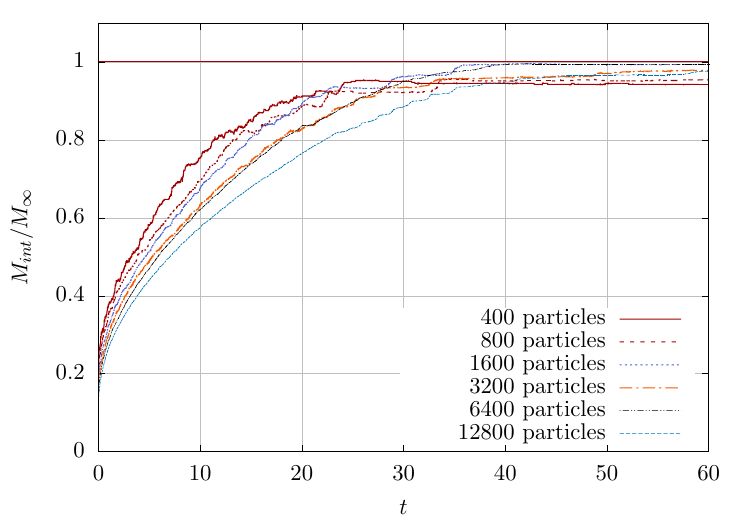}
    \caption{\footnotesize
    Parallelepiped example.
    Time evolution of the mass, relative to the desired value $M_\infty$ for various ideal numbers of particles $n_\infty$.}
    \label{fig:fig8}
\end{figure}
\begin{figure}[h]
    \centering
    \includegraphics[width=0.7\textwidth]{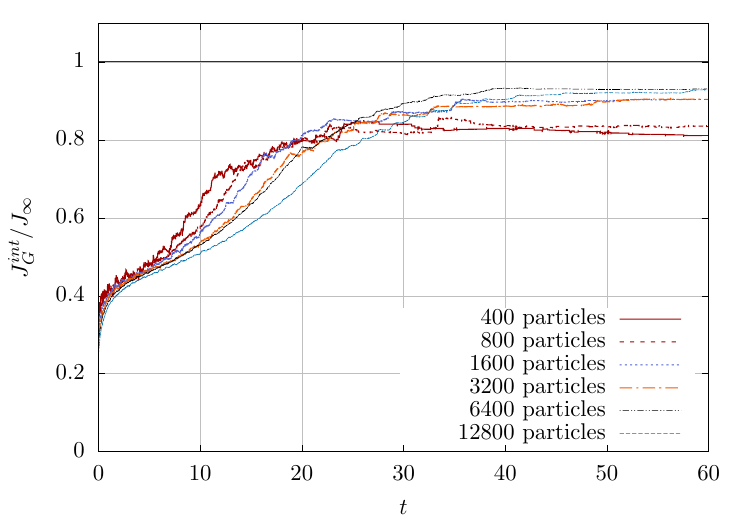}
    \caption{\footnotesize
    Parallelepiped example.
    Time evolution of the polar inertia, relative to the desired value $J_\infty$ for various ideal numbers of particles $n_\infty$.}
    \label{fig:fig9}
\end{figure}

We set the diffusivity $\kappa=1$ and the total analysis time $T=60$. The Dirichlet boundary condition is implemented using a boundary layer of size $2b=2\sqrt{L\,\Delta x}$ with
\begin{equation}
    \Delta x = \left( \frac{3\, m_p}{4\,\pi\,\rho_\infty} \right)^{1/3} \, .
\end{equation}
For the particle Gaussian profile, we adopt $\beta=2/\Delta x^2$. Consequently, the number of particles in the steady state regime $n_\infty$ is approximately $M_\infty/m_p$.

\begin{figure}[t]
    \centering
    \includegraphics[width=0.7\textwidth]{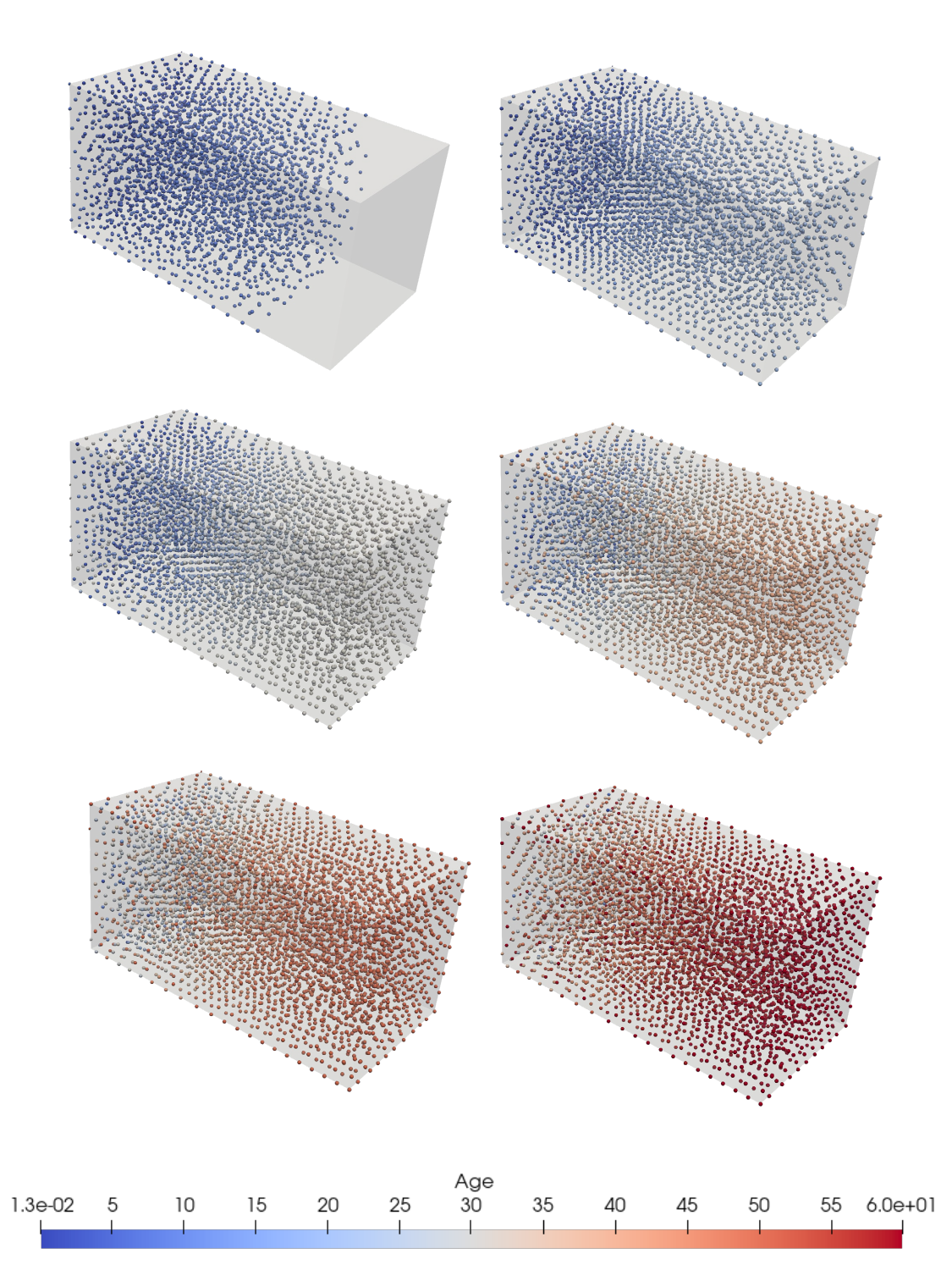}    
    \caption{\footnotesize Parallepiped example.
      Snapshots of the filling process ($n_p=3200$). From top to bottom, left to right, $t=10\cdot k, k=1,2,\ldots6$. Colors refer to particle age since
      insertion.}
    \label{fig-longbrick-sn}
\end{figure}

\begin{figure}[h]
    \includegraphics[width=0.7\textwidth]{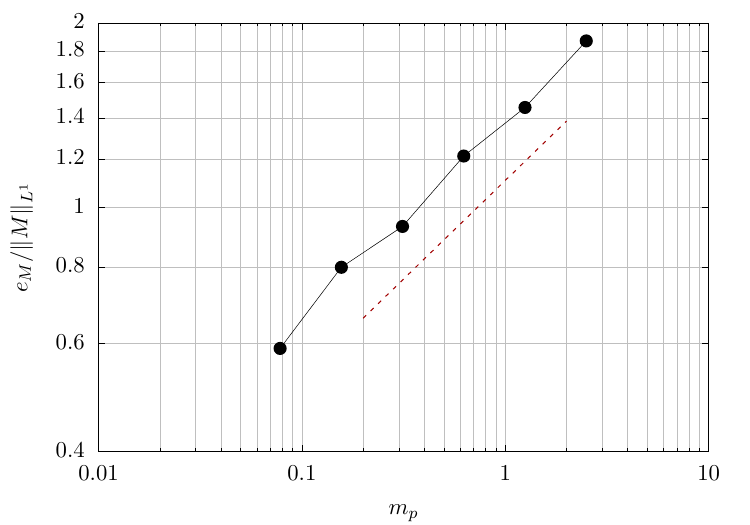}
    \caption{\footnotesize
    Parallelepiped example. Relative error in the $L^1$ norm of the total mass. The Broken line represents the best-fitting power law ($e_{M}\approx m_p^{0.32}$). Normalization factor is the fitted asymptotic value $\|M\|_{L^1(0,T)}=43060$.} 
    \label{fig:fig11}
\end{figure}

\begin{figure}[h]
    \centering
    \includegraphics[width=0.7\textwidth]{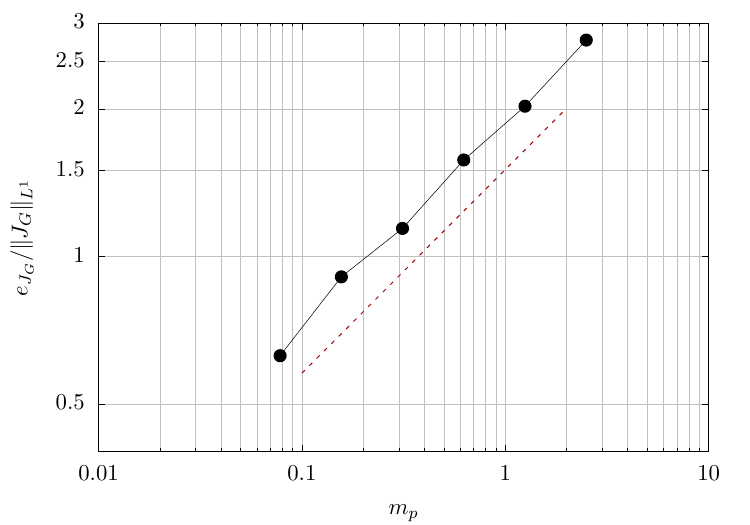}
    \caption{\footnotesize
     Parallelepiped example. Relative error in the $L^1$ norm of the total polar inertia. The broken line represents the best-fitting power law ($e_{J_G}\approx m_p^{0.42}$). Normalization factor is the fitted asymptotic value $\|J_G\|_{L^1(0,T)}=27660$).}
    \label{fig:fig12}
\end{figure}

\begin{table}[ht]
    \centering
    \begin{tabular}{rcccccrr}
    \hline
    $n_\infty$ & $m_p$ & $\Delta x$ & $\beta$ & $n_{\rm t}$ & $n_{\rm t}/n_{\infty}$ \\
     \hline
      \hline
   400 &      2.5 &  0.10607847 & 177.73604 &   908 &      2.27 \\
   800 &     1.25 & 0.084194539 & 282.13837 &  1593 &   1.99125 \\
  1600 &    0.625 & 0.066825250 & 447.86674 &  2879 &  1.799375 \\
  3200 &   0.3125 & 0.053039236 & 710.94414 &  5154 &  1.610625 \\
  6400 &  0.15625 & 0.042097269 & 1128.5535 &  9635 & 1.5054688 \\
 12800 & 0.078125 & 0.033412625 & 1791.4670 & 17579 & 1.3733594 \\
      \hline
    \end{tabular}
    \vskip 5pt
    \caption{\footnotesize
    Numerical parameters and results for the parallelepiped example.
    $n_{\rm t}$ is the total number of particles, including the
    external layer, in the stationary regime, $n_{\rm t}/n_\infty$ is
    the particle overhead due to the Dirichlet condition on one
    of the faces of the domain.}
    \label{table-parall}
\end{table}
We consider six sets of particles with $n_\infty$ ranging between 400 and 12800, see Table~\ref{table-parall}. As in the previous example, the boundary layer where the Dirichlet conditions are imposed uses an external buffer of particles, which adds some overhead to the computational cost of the simulation. However, as the number of particles increases, the relative extra cost reduces and eventually becomes negligible.

Fig.~\ref{fig:fig8} shows the time evolution of the internal mass relative to the desired value $M_\infty$ for the six cases described in Table~\ref{table-parall}. These results demonstrate that the algorithm successfully enforces the desired density for all particle distributions, even when mixed boundary conditions are applied. Finally, for the six cases, Fig.~\ref{fig:fig9} illustrates the evolution of the polar inertia of the particle distribution inside the parallelepiped relative to the theoretical value
\begin{equation}
    J_\infty = \frac{1}{12} M ( L^2 + 2 B^2) = 500 \ .
\end{equation}

Fig.~\ref{fig-longbrick-sn} shows six snapshots that describe the filling process of the parallelepiped. The particles are colored according to their age, defined by the time since injection. In particular, `older' particles are always positioned away from the boundary layer, located to the left of each box in the figure.

The convergence of the method with the Dirichlet and zero-flux Neumann boundary conditions is illustrated in Fig.~\ref{fig:fig11}, where the error in the $L^1$ norm of the total mass, as defined in Eq.~\eqref{eq-l1-sphere1}, is shown. Similarly, Fig.~\ref{fig:fig12} depicts the error in the $L^1$ norm of the polar inertia of the mass distribution relative to the origin, as defined in Eq.~\eqref{eq-l1-sphere2}.

\subsection{Flow in a circular pipe with asymmetric inlet and outlet}
\label{ssec:pipe}

\begin{figure}[h]
    \centering
    \includegraphics[width=0.8\textwidth]{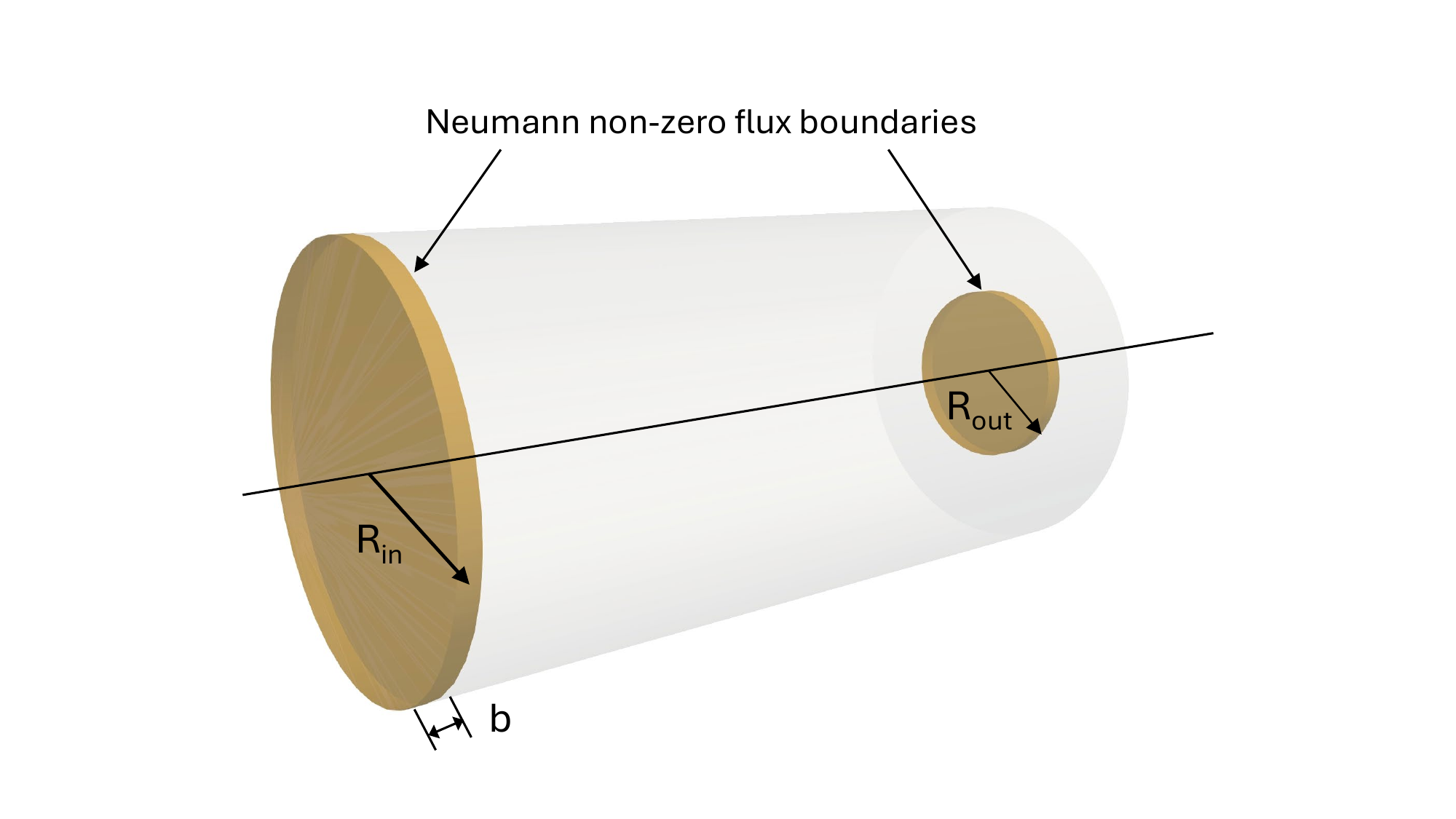}
    \caption{\footnotesize
    Circular prism domain used for the verification of the particle method in a problem with assigned inlet and outlet fluxes. Inlet flux on the left, outlet flux on the right.}
    \label{fig:fig15}
\end{figure}
We conclude with an example to assess Neumann boundary conditions. We consider a cylinder of length $L=2$ and radius $R_{\rm in}=0.5$, see Fig.~\ref{fig:fig15}. The entire left face of the cylinder is an inlet surface, with constant flux assigned $f_{\rm in}=5000$. On the right face, a circular outlet surface of radius $R_{\rm out}=0.25$ is placed at the center with assigned constant flux $f_{\rm out} = f_{\rm in} R^2_{\rm in}/R^2_{\rm out}=20000$, so that the influx and outflux are equal. We assign the diffusivity $\kappa = 1$ and the total time $T=6$.

The domain is initially filled with $n_{\rm initial}$ particles of total mass $M=1000$, providing an initial average density $\rho_0 = M/|\Omega|$. Given
\begin{equation}
    \Delta x = \left( \frac{3}{4} \, \frac{R^2_{\rm in}\,L}{n_{\rm initial}} \right)^{1/3} 
\end{equation}
we set the boundary layer thickness $b = \sqrt{R\,\Delta x}$, the width parameter $\beta = 1/\Delta x^2$.

At the two inlet and outlet boundaries, the number of particles to be injected or removed is
\begin{equation}
    n_{\rm in} =  \left\lfloor f_{\rm in} \frac{\pi R^2_{\rm in} \Delta t}{m_p}\right\rfloor \, ,
    \qquad
    n_{\rm out} =  \left\lfloor f_{\rm out} \frac{\pi R^2_{\rm out} \Delta t}{m_p} \right\rfloor\,.
\end{equation}
which, for the choice of the fluxes, are identical. Particles are injected into the boundary layer or removed from the boundary layer at random positions. 
In the simulations we consider an increasing number of initial particles ranging between 1000 and 32000, see Table~\ref{table-pipe}.

\begin{table}[]
    \centering
    \begin{tabular}{rcccccrr}
        \hline
         $n_{\rm initial}$  & $m_p$ & $dx$ & $b$ & $dt$ & $\beta$ & $n_{\rm bl}$ & $n_{\rm final}$ \\
         \hline
         \hline
         1000 & 1. & 0.0721 & 0.190 & 0.005200 & 192.30 & 20 & 1347 \\
         2000 & 0.5 & 0.0572 & 0.169 & 0.003275 & 305.26 & 25 & 2365 \\
         4000 & 0.25 & 0.0454 & 0.150 & 0.002063 & 484.57 & 32 & 4240 \\
         8000 & 0.125 & 0.0360 & 0.134 & 0.001300 & 769.20 & 40 & 7617 \\
         16000 & 0.0625 & 0.0286 & 0.120 & 0.000819 & 1221.03 & 51 & 14005 \\
         32000 & 0.03125 & 0.0227 & 0.107 & 0.000519 & 1938.26 & 64 & 25835 \\
         \hline
    \end{tabular}
    \vskip 5pt
    \caption{\footnotesize
    Numerical parameters and results for the pipe example with assigned fluxes. Here $b = \sqrt{L \Delta x}$ and $\beta = \Delta x^{-2}$. $n_{\rm initial}$ is the initial number of particles necessary to achieve a prescribed uniform density $\rho_0$ in the domain. $n_{\rm bl}$ is the number of particles that are injected into or removed from the two inlet and outlet boundary layers. $n_{\rm final}$ is the number of particles obtained at the steady state.}
    \label{table-pipe}
\end{table}

\begin{figure}[t]
    \centering
    \includegraphics[width=0.7\textwidth]{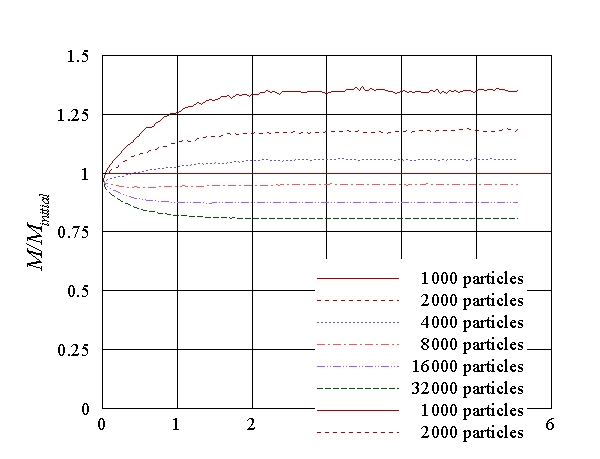}
    \caption{\footnotesize
    Pipe example, with inlet and outlet fluxes, using $b = \sqrt{R\, \Delta x}$ and $\beta = 1/\Delta x^2$. Time evolution of the mass, normalized by the initial mass, always set $M=1000$, for various initial numbers of particles $n_{\rm initial}$.}
    \label{fig:fig16}
\end{figure}

Fig.~\ref{fig:fig16} shows the time history of the total mass of the system, normalized with the initial mass, for the six cases, showing that the total mass in steady state is an outcome of the analysis.
\begin{figure}[t]
    \centering
    \includegraphics[width=0.7\textwidth]{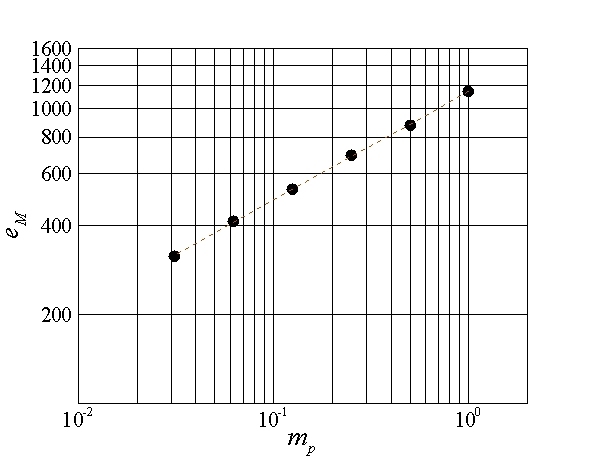}
    \caption{\footnotesize
    Pipe example, with inlet and outlet fluxes, using $b = \sqrt{R\, \Delta x}$ and $\beta = 1/\Delta x^2$. $L^1$-norm error of the total mass at the steady state. Circles denote the values $L^1$-norm of the total mass computed at the steady state. The broken line represents the best-fitting power law showing convergence ($e_M\approx m_p^{0.37}$). Normalization factor is the fitted asymptotic value $\|M\|_{L^1(0,T)}=619$. }
    \label{fig:fig17}
\end{figure}

\begin{figure}[t]
    \centering
    \subfigure[$n_{\rm initial}$ = 2000]{\includegraphics[width=0.7\textwidth]{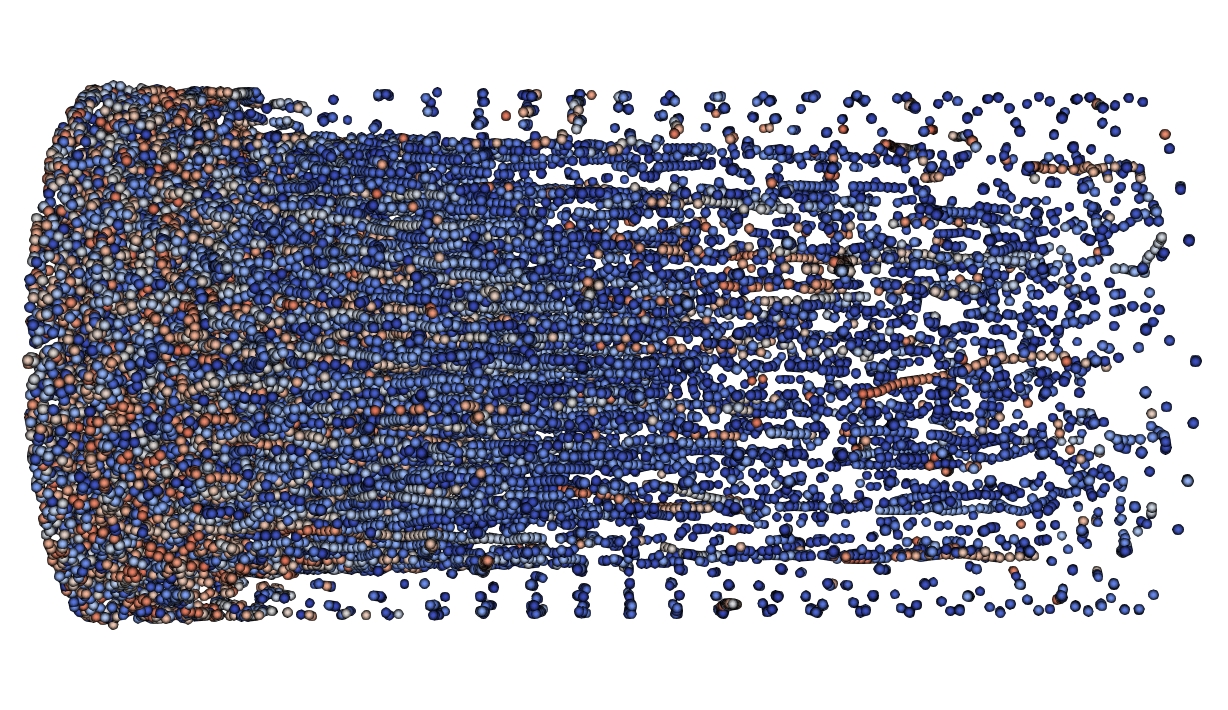}}
    \subfigure[$n_{\rm initiall}$ = 32000]{\includegraphics[width=0.7\textwidth]{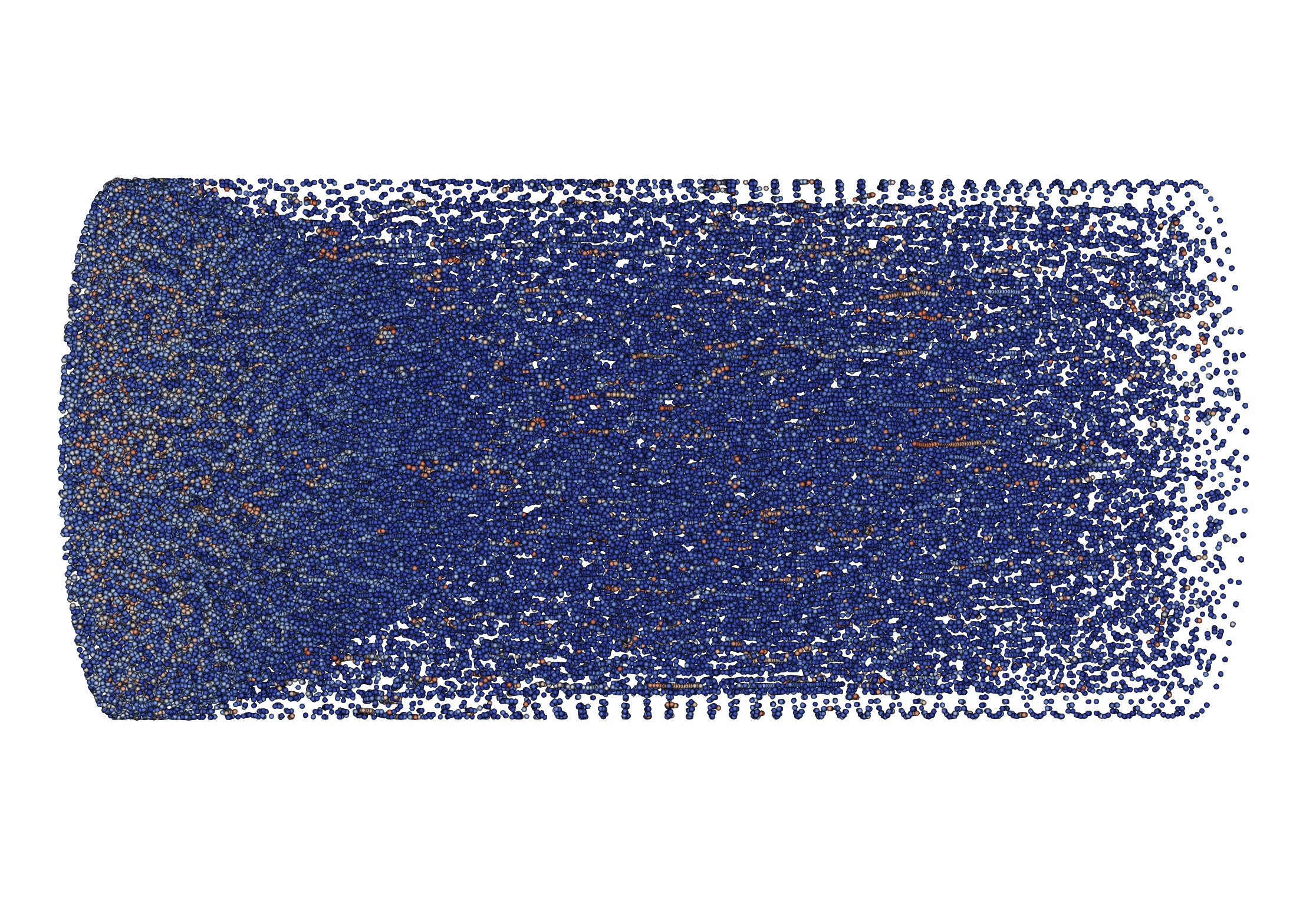}}
    \caption{\footnotesize
    Pipe example. Cases $n_{\rm initial}=2000$ and $n_{\rm initial} = 32000$.  Visualization of the sequential positions occupied by the particles at time 6.}
    \label{fig:fig18}
\end{figure}

Fig.~\ref{fig:fig18} visualizes the sequence of the positions occupied by the particles at the final time of the analysis, showing the average path followed by the particles in the two cases corresponding to the initial particles 2000 and 32000.

\section{Summary and concluding remarks}
\label{sec:conclusions}

We have extended the finite-particle method of mass transport of \cite{Pandolfi:2023} to account for general mixed boundary conditions. In addition to the geometrically-exact treatment of advection coupled to diffusive dynamics characteristic of particle methods, the extension to general boundary conditions enables consideration of commonly encountered effects such as equilibrium with the environment and inlet/outlet conditions. General boundary conditions are enforced by introducing an adsorption/depletion boundary layer wherein particles are added or removed in accordance with the boundary conditions. We showcase the range and scope of the method through a number of examples of application, including absorption of particles into a sphere and pipe flow through pipes of square and circular cross sections, with and without occlusions.

The fundamental tenet of the present approach--and of all particle methods generally--is that mass diffusion is most naturally understood as a problem of transport of {\sl measures}, instead of the traditional view of diffusion as a problem of evolution in linear spaces of functions. That means that the solution is represented by a measure, i.~e., a rule for computing spatial averages of continuous functions representing quantities of interest. Consistent with this view, particle mobility is represented, in discrete time, by the Wasserstein distance between consecutive configurations of the particle ensemble, in the spirit of Wasserstein gradient flows \cite{JordanKinderlehrerOtto1997}. Likewise, entropic forces are formulated by an appeal to the Kullback-Leibler functional, which is natural for regular measures, following mollification of the particles.

This regularization of the entropy is closely related to the `blob' method of Carrillo {\sl et al.} \cite{carrillo2017, carrillo2019} and has been analyzed in detail by Conti {\sl et al.} \cite{Conti:2023}. In the present work, we investigate convergence numerically by way of example. As expected from the measure-theoretical framework on which the method is predicated, in all cases, the numerical solution is observed to converge weakly, or in the sense of local averages.

We conclude by remarking on the disarming implementation simplicity of the method, as attested to by the computer programs provided in the supplementary materials. However, much of the simplicity of the implementation rests on a na\"ive evaluation of the diffusive forces (\ref{j4K0Ht}), an $O(n^2)$ operation that, for large $n$, becomes the main computational bottleneck. Calculations can be greatly accelerated by noting that, due to the Gaussian decay of the particle profiles, the sums can be restricted to local neighborhoods of size $\ell=1/\sqrt{\beta}$ without appreciable loss of
accuracy. An overview of techniques for solving the local-neighborhood problem can be found in \cite{tapiafernandez2017ip, tapiafernandez2021ba}, and a particular implementation based on spatial decomposition of the analysis domain into regular cells is provided in \cite{Pandolfi:2023}.

\section*{Supplementary material}
The code used for the examples in this article is freely available in the public repository 
\\
\centerline{\texttt{https://gitlab.com/ignacio.romero/finite-particles}}

\section*{Acknowledgements}
AP is grateful for the support of the Italian National Group of Physics-Mathematics (GNFM) of the Italian National Institution of High Mathematics ``Francesco Severi'' (INDAM). IR acknowledges the support of the Spanish Ministry of Science and Innovation under project PID2021-128812OB-I00. MO gratefully acknowledges the support of the Deutsche Forschungsgemeinschaft (DFG, German Research Foundation) {\sl via} project 211504053 - SFB 1060; project 441211072 - SPP 2256; and project 390685813 -  GZ 2047/1 - HCM.

%
%
%
%
%


\end{document}